\documentclass[12pt, leqno]{article}
\usepackage{hyperref}
\usepackage{amsmath}
\usepackage{amsfonts}
\usepackage{amssymb}
\usepackage{ifpdf}
\usepackage{graphicx}
\usepackage{color}
\usepackage{color}

\newtheorem{thm}{Theorem}[section]

\newtheorem{lem}{Lemma}[section]

\numberwithin{equation}{section}

\newcommand{\dx}{\,{\rm d}x}

\newcommand{\RR}{\mathbb{R}}
\newcommand{\re}{\mathbb{R}}
\newcommand{\ren}{\mathbb{R}^N}

\newcommand{\la}{\lambda}
\newcommand{\ve}{\varepsilon}
\newcommand{\vp}{\varphi}
\newcommand{\sign}{\mbox{sign\,}}

\newcommand{\dint}{\displaystyle\int}
\newcommand{\diint}{\displaystyle\iint}
\renewcommand{\dfrac}{\displaystyle\frac}
\def\qed{\,\unskip\kern 6pt \penalty 500
\raise -2pt\hbox{\vrule \vbox to8pt{\hrule width 6pt
\vfill\hrule}\vrule}\par}
\definecolor{darkblue}{rgb}{0.05, .05, .65}
\definecolor{darkgreen}{rgb}{0.05, .70, .05}
\definecolor{darkred}{rgb}{0.8,0,0}
\def\qed{\unskip\kern 6pt \penalty 500
\raise -2pt\hbox{\vrule \vbox to8pt{\hrule width 6pt
\vfill\hrule}\vrule}\par}
\begin{document}
\title{\textbf{Regularity of solutions of the \\  fractional porous medium flow \\
with exponent $1/2$}}
\author{\Large Luis Caffarelli \footnote{University of Texas; \sl caffarel@math.utexas.edu}
 \\\Large Juan Luis Vazquez\footnote{Universidad Aut\'onoma de Madrid; \sl juanluis.vazquez@uam.es}\\}

%\date{Draft version, \today}
\maketitle

\begin{abstract}
We study the regularity of a porous medium equation with nonlocal diffusion effects
given by an inverse fractional Laplacian operator. The precise model  is \
$u_t=\nabla\cdot(u\nabla (-\Delta)^{-1/2}u).$
For definiteness, the problem is posed in $\{x\in\ren, t\in \re\}$ with nonnegative initial data $u(x,0)$ that are integrable and decay at infinity. Previous papers have established the existence of mass-preserving, nonnegative weak solutions satisfying energy estimates and finite propagation, as well as the boundedness of nonnegative solutions with $L^1$ data, for the more general family of equations $u_t=\nabla\cdot(u\nabla (-\Delta)^{-s}u)$, $0<s<1$.

Here we establish the  $C^\alpha$ regularity of such weak solutions in the difficult fractional exponent case $s=1/2$.
For the other fractional exponents $s\in (0,1)$ this H\"older regularity has been proved in \cite{CSV}. The method combines delicate De Giorgi type estimates with iterated geometric corrections that are needed to avoid the divergence of some essential energy integrals  due to fractional long-range effects.
\end{abstract}

\vfill


{\bf Keywords and phases.} Porous medium equation, fractional Laplacian, nonlocal diffusion operator, H\"older regularity.

\medskip

\newpage

%%%%%%%%%%%%%%%%%%%%%%%%%%%%%%%%%%%%%%%%%%%%%%%%%%%%%%%%%%%%%%%%%%
\section{Introduction}
\label{sec.intro}

This paper is devoted to study the regularity properties of nonnegative weak solutions of a model of porous medium equation that includes nonlocal effects through an integral operator of fractional Laplacian type, thus allowing to account for long-range effects. In  \cite{CaVa09}  we have introduced the following concrete model
\begin{equation}\label{eqs}
u_t=\nabla\cdot(u\nabla p), \quad p =(-\Delta)^{-s}u, \ 0<s<1\,,
\end{equation}
posed for $x\in \RR^N$, $N\ge 2$, and $t>0$, with  initial conditions
\begin{equation}\label{eq.ic}
u(x,0)=u_0(x), \quad x\in \RR^N,
\end{equation}
where $u_0$ is a nonnegative and integrable function in $\RR^N$ decaying as $|x|\to\infty$.
We will refer to Equation \eqref{eqs} as the FPME (for fractional porous medium equation), but we remind the reader that another natural candidate to the denomination of fractional porous medium equation has been recently introduced and studied in the literature, see \cite{PQRV1, PQRV2}, and both models have quite different properties, cf. \cite{JLVsurvey2}.

We have proved existence of nonnegative weak solutions for the Cauchy problem \eqref{eqs}-\eqref{eq.ic}, enjoying a number of additional properties. The long-time  behaviour of such solutions was described in a second paper \cite{CaVa11} using an entropy method that leads to a fractional obstacle problem for the asymptotic profile. In a more recent paper with F. Soria \cite{CSV} we have addressed the questions of boundedness and H\"older regularity of such solutions. Boundedness was established for all exponents $0<s<1$ with a clean estimate of the form
\begin{equation}\label{form:L-inf}
\sup_{x\in\mathbb{R}^n}|u(x,t)|\le C\,t^{-\alpha }\|u_0\|_{L^1(\mathbb{R}^n)}^{\gamma}
\end{equation}
with precise exponents $\alpha=N/(N+2-2s)$, $\gamma=(2-2s)/(N+2-2s)$ (as corresponds to dimensional analysis), and a constant $C$ that depends only on $N$ and $s\in (0,1)$. For $s=1/2$ the exponents are $\alpha=N/(N+1)$ and $\gamma=1/(N+1)$.

Concerning regularity, we found in \cite{CSV} a $C^\alpha$ modulus of continuity for all exponents $s\in (0,1)$, $s\ne 1/2$, for bounded, nonnegative, and integrable weak energy solutions defined in a strip of space-time.

The proof of the $C^\alpha$ result in the range $0<s<1/2$ uses a number of techniques
that are becoming classical in the study of regularity of nonlocal diffusion problems, but the process is long and delicate since we must take into account  the nonlinearity with possible degeneracy, as well a the long range interaction carried by the kernel of the fractional operator.

The regularity result in the case $1/2<s<1$  is more difficult due to the last mentioned effect and the method proposed in \cite{CSV} uses a geometrical transformation to absorb the uncontrolled growth of one of the integrals that appear in the iterated energy estimates on which the regularity method is based. In other words, we control the possible divergences using transport.

\medskip

\noindent {\bf The case $s=1/2$.} It appeared as a borderline case for those methods and was left out in previous papers since it needs a new version of the technique that takes some space and involves quite careful iterative constructions. It is the purpose of this paper to perform the proof of H\"older regularity in such critical case in full detail. Therefore, in the sequel  we will concentrate on the Cauchy problem for equation \eqref{eqs} with exponent $s=1/2$. We may write the equation as
\begin{equation}\label{eq1}
u_t+\nabla\cdot(u \,{\bf v})=0, \quad {\bf v}=-\nabla  (-\Delta)^{-1/2}u,
\end{equation}
posed for $x\in \RR^N$, $N\ge 1$, and $t>0$, with  initial condition \eqref{eq.ic}. Since we are interested in regularity, and in view of estimate \eqref{form:L-inf},  we assume without loss of regularity that $u$ is bounded.
Moreover, the velocity field can be written in the present case as ${\bf v}=-{\bf R}u$, where $\bf R$ is the Riesz transform with components $R_k=\partial_k (-\Delta)^{-1/2}$  (and symbol $i \xi_k/|\xi|$).
In terms of singular integrals we may write
\begin{equation}
{\bf v}_k(x,t)=  \, c_N \,\text{P.V.}\int \frac{(x_k-y_k)\,u(y,t)}{|x-y|^{N+1}}\,dy\,,
\qquad k=1,2,\cdots, N,
\end{equation}
with $c_N= 1/(\pi \omega_{N-1})=\Gamma((N+1)/2)\,\pi^{-(N+1)/2} $.  Here is our main result:

\begin{thm}\label{mainthm} Let $u\ge 0$ be a bounded weak energy solution of equation \eqref{eq1} with  initial data $u_0$ that is a nonnegative and integrable function in $\RR^N$ decaying exponentially as $|x|\to\infty$, and assume that $u$  defined in a space-time domain $S=\RR^N\times [0,T]$. Then  $u$ is $C^\alpha$ continuous in any subdomain $Q\subset \RR^N\times [\tau,T]$ with some H\"older exponent \ $\alpha(N,s)\in (0,1)$  and a constant that depends also on the bound for $u_0$ in $L^1(\ren)$ and the minimal time $\tau$ of \ $Q$.
\end{thm}

This completes the $C^\alpha$ regularity result for all $s\in (0,1)$.
The main novelty is as follows: in the most delicate situation, the degenerate case where $u$ approaches zero, we perform  a careful iteration analysis that combines consecutive applications of  scaling and a geometrical transformation. The end result of  H\"older regularity is derived from the estimate of the size of the solution oscillation in a sequence of shrinking and distorted cylinders.

It is interesting to recall that the fractional Laplacian exponent $s=1/2$ is the most natural case in the whole family $s\in (0,1)$, both for the number and interest of the applications
%%% ojo add here
in different scientific contexts, and  for the simplicity of the definition in the sense of Caffarelli-Silvestre extension, \cite{CS07}. However, this is precisely the case where we encounter the biggest difficulties in the regularity theory. This is maybe due to the fact that the net order of differentiation of the diffusive term, $\nabla\cdot(u\nabla (-\Delta)^{-1/2}u), $ is just 1, which makes it formally the same order as the convective equations, see \eqref{eq1}.
The conclusion of this paper is, roughly speaking, that the diffusive character is still dominant.

It is interesting to compare  the present paper with the paper \cite{Caffarelli-Vasseur} by Vasseur and one of the authors where
boundedness and regularity is established for the geostrophic model
\begin{equation}
\partial_t u + {\bf v}\cdot \nabla u+ (-\Delta)^{1/2} u=0\,,
\end{equation}
with $\nabla\cdot {\bf v}= 0$. There are two main differences: the convective velocity is not necessarily divergence-free in our case, and the last term in our model is nonlinear and degenerate
since we can write our equation as \ $\partial_t u + {\bf v}\cdot \nabla u+ u(-\Delta)^{1/2} u=0$. This explains maybe the very involved analysis that we will have to perform.

We refer to the above-mentioned papers and  \cite{JLVAbel} for motivation and background on fractional porous medium equations. For basic information on fractional Laplacians see e.g. \cite{Landkof, Stein70, Valdinoci}. On the other hand, the case $N=1$ of our equation is a bit special. The equation can be written as $u_t+(H(u)u)_x=0$ where $H(u)$ denotes the Hilbert transform of $u$. This simpler case was treated in \cite{Head72, dtdcnm04, BKM} in a problem coming from dislocation theory, and in \cite{BLM96, CC08, CCCF} as a simplified model in fluid mechanics.  Our result applies to that 1D case with minor modifications in the proofs, see further comments at the end of the paper.

\medskip

\noindent {\bf Organization.}  We collect the necessary preliminaries in Section \ref{sec.prelim} and present the strategy of proof in Section \ref{sec.holder}, where the similarities and differences with the cases covered in \cite{CSV} are explained. Putting this into effect, Sections \ref{sec51} to \ref{sec.3lemma} contain the main lemmas that will be needed in versions adapted to the case $s=1/2$, and Section \ref{sec.iterproc} reviews the iterative procedure to regularity already introduced in \cite{CSV}. This is done with the detail needed  to prepare the reader for the modifications that are needed to avoid divergent integrals. After the presentation of the transport correction method used for $s>1/2$  (Section \ref{s-large}), the main contribution of the paper is presented in Section \ref{s=12}: it is an iterative method with geometrical transport corrections which allows to obtain convergent and controlled integrals for the regularity estimates near a degenerate point (i.\,e., where the solution vanishes in some sense); it involves successive corrections based on a sequence of cutoffs, plus summation of all the corrections via a geometrical series estimate. Three short last sections add different complements: one for the non-degenerate case, which is easier; one for the equation in 1D, which has some peculiarities, and one for general data, that recalls the extension of a result proved in \cite{CSV}.

\medskip

\noindent {\bf Notations.}  As a general rule, we will use the notations of paper \cite{CSV}. Thus, $(-\Delta)^{s}$ with $0<s<1$ denotes the fractional powers of the Laplace operator defined on the Schwartz class of functions in $\RR^N$ by Fourier transform, and extended in a natural way to functions in the Sobolev space $H^{2s}(\RR^N)$.   The inverse operator is denoted by ${\cal L}_s=(-\Delta)^{-s}$ and can be realized by convolution
\begin{equation*}
{\cal L}_s u=L_s\star u, \qquad L_s(x)=c(N,s)|x|^{-N+2s}.
\end{equation*}
${\cal L}_s$ is a positive self-adjoint operator. We will write ${\cal H}_s={\cal L}_s^{1/2}$ which has kernel $L_{s/2}$.  For $s=1/2$ we simply write ${\cal L}=(-\Delta)^{-1/2}$, ${\cal H}=(-\Delta)^{1/4}$. For functions that depend on $x$ and $t$, convolution is applied for every fixed $t$ with respect to the space variables. We then use the abbreviated notation $u(t)=u(\cdot,t)$.

For a measurable  $u\ge 0$ and  for $k>0$ we denote by $u_k^+=(u-k)_+=\max\{u-k,0\}$, and $u_k^-=\min\{u-k,0\}$ in such a way that $u_k^+\ge 0\ge u_k^-$, the supports of $u_k^+$ and $u_k^-$ agree only on points where $u=0$, and also $u=k+u_k^+ + u_k^-$. We will use similar notations: $u_\varphi^+=(u-\varphi)^+$, $u_\varphi^-=(u-\varphi)^-$ when $\varphi$ is a function and not just a constant, and then we may split $u$ as follows: $u={\varphi}+u_{\varphi}^+ + u_{\varphi}^-$. The notation $|\Omega|$  means the Lebesgue measure of the set $\Omega$.

%%%%%%%%%%%%%%%%%%%%%%%%%%%%%%%%%%%%%%%%%%%%%%%%%%%%%%%%%%%%%%%%%%
\section{Preliminaries. Existence and basic estimates}
\label{sec.prelim}

From now on we will concentrate on the Cauchy problem with exponent $s=1/2$.
We may write the equation in the form \eqref{eq1}. We will also write equation \eqref{eq1} in the form
\begin{equation}\label{eq1.p}
u_t=\nabla\cdot(u \,\nabla p), \quad p={\cal L }u := (-\Delta)^{-1/2}u,
\end{equation}
and call $p$ the pressure of the solution. The two forms are clearly equivalent. We work on dimensions $N\ge 2$
and extend the results to dimension $N=1$ in Section \ref{sect.1d}.

\medskip

\noindent {\bf Concept of solution.} We will work in the class of weak  nonnegative solutions that have some additional properties. In the present situation we start with a slightly polished  version of the definition of \cite{CaVa09} since our solutions will always be bounded.

\noindent $\bullet$ By a  {\sl weak solution} of Problem {\rm (\ref{eq1})-\eqref{eq.ic}} defined in  a space-time domain \ $Q_T=\RR^N\times (0,T)$ \ with initial data \ $u_0\in L^1(\RR^N)\cap L^\infty(\RR^N)$ we understand a nonnegative function $u\in L^1(Q_T)\cap  L^\infty(Q_T)$, such that ${\cal L}u=(-\Delta)^{-1/2} (u)\in L^1(Q_T)$, and  the identity
\begin{equation}
\iint u\,(\eta_t-\nabla {\cal L}(u)\cdot\nabla\eta)\,dxdt+ \int
u_0(x)\,\eta(x,0)\,dx=0
\end{equation}
holds for all  continuously differentiable test functions $\eta$ in $Q_T$  that are compactly supported in the space variable and vanish near $t=T$.

\medskip

\noindent {\bf Existence of solutions.} The following results have been proved in \cite{CaVa09}:
For any  $u_0\in  L^\infty(\ren)\cap L^\infty(\ren)$, $u_0\ge 0$, and such that
\begin{equation}
u_0(x)\le A\,e^{-a|x|} \qquad \mbox{for some $A,a>0$}\,,
\end{equation}
 there exists a weak solution $u\ge 0$ of Equation \eqref{eq1} with initial data $u_0$.  Moreover, for all $t>0$ we have the conservation of mass:
\begin{equation}
\int_{\ren} u(x,t)\,dx=\int_{\ren} u_0(x)\,dx,
\end{equation}
as well as the $L^\infty$ bound: $\|u(t)\|_\infty \le \|u_0\|_\infty$. The constructed solution decays exponentially as $|x|\to\infty$. According to the smoothing effect formula \eqref{form:L-inf}, it decays also in time like $O(t^{-N/(N+1)})$ for large times.

\medskip

\noindent {\bf Other properties of the constructed solutions.} Here are some of the most useful
properties that are known.

\noindent $\bullet$ {Translation invariance.} The equation is invariant under translations in space and time, and this property reflects on the set of weak solutions.

\medskip

\noindent $\bullet$ {Scaling.} Moreover, the equation is invariant under a subgroup of the group of dilations in $(u,x,t)$, and this implies a scaling property for the set of solutions. To be precise, if $u(x,t)$ is a weak solution as described in the existence theorem, with initial data $u_0(x)$, and $A,B,C$ are positive constants, then $\widehat u(x,t)= A\, u(Bx, Ct)$ is again a weak solution on the condition that $C=AB$. It has initial data $\widehat u_0(x)=Au_0(Bx)$.

\medskip

\noindent $\bullet$ {Finite propagation:}
Compactly supported initial data $u_0(x)$ give rise to solutions $u(x,t)$ that have the same
property for all positive times, i.e., the support of $u(\cdot,t)$
is contained in a finite ball $B_{R(t)}(0)$ for any $t>0$.

 \medskip

\noindent  $\bullet$ A standard comparison result for parabolic equations does not work in general. This is one of the main technical difficulties in the study of the equation.

\medskip

\noindent  {\bf Energy Properties. }  The boundedness and regularity analysis performed in \cite{CSV} uses in an essential way regularity properties that go beyond the definition of weak solution, but they are satisfied  by the solutions constructed in \cite{CaVa09}. We call them energy properties since they involve bilinear forms and integrals of (fractional) derivatives.

\noindent  $\bullet$  The first energy inequality holds in the form
\begin{equation}
\dint_0^t\dint_{\ren} |\nabla (-\Delta)^{-1/4} u|^2\,dxdt +  \dint_{\ren}  u(t)\log(u(t))\,dx
 \le \dint_{\ren} u_0\log(u_0)\,dx\,,
\end{equation}
This estimate allows to control (uniformly in $t$) the norm of $u$ in ${W}^{1/2,2}$ and implies compactness in space for the solutions.

\noindent  $\bullet$ The second energy estimate of \cite{CaVa09} says that for all $0<t_1<t_2<\infty$
\begin{equation}
\int_{t_1}^{t_2}\int_{\ren} u\,|  {\bf R}u|^2\,dxdt+ \frac 12\int_{\ren} |(-\Delta)^{-1/4} u(t_2)|^2\,dx\le
\frac 12\int_{\ren} |(-\Delta)^{-1/4}(u(t_1)|^2\,dx\,.
\end{equation}
This estimate is not so important in a context in which solutions are integrable and bounded. More important is the following observation.

\medskip

\noindent  $\bullet$ { \bf General energy property. Bilinear form.} (i) The proofs of the basic lemmas will use energy inequalities of the following form: for any $F$ smooth and such that $f=F'$ is bounded and nonnegative, we have   for every $0\le t_1\le t_2\le T$
\begin{equation*}\label{entro}
\begin{array}{ll}
\dint F(u(t_2))\, dx -\dint F(u(t_1))\, dx & = - \dint_{t_1}^{t_2}\dint \nabla [f(u)] u \nabla p\, \,dx\, dt= \\
& -\dint_{t_1}^{t_2}\int \nabla h(u) \nabla (-\Delta)^{-1/2} u\, dx\, dt\,,
\end{array}
\end{equation*}
where $p=(-\Delta)^{-1/2}u$ and  $h$ is a function  satisfying  $h'(u)= u\,f'(u)$. A natural candidate for $h(u)$ is $u$ or a truncation of it, so that $f(u)$ must be a logarithm or a variant of it, as we will see below. We can write the last integral as a bilinear form
\begin{equation}
\int \nabla h(u) \nabla (-\Delta)^{-1/2} u\, dx = \mathcal{B} (h(u), u)
\end{equation}
and this bilinear form $\mathcal{B}$ is defined on the  Sobolev space $W^{1,2}(\ren)$ by
\begin{equation}\label{B}
\mathcal{B}(v,w)= C_{N} \iint  \nabla v(x)\frac{1}{|x-y|^{N-1}} \nabla w(y)\, dx \, dy,
\end{equation}
where ${\cal N}(x,y)= C_{N}|x-y|^{-(N-1)}$ is the kernel of operator $(-\Delta)^{-1/2}$.
After some well-known transformations, we also have
\begin{equation}\label{BB}
\mathcal{B} (v,w)= C'_{N} \iint(v(x)-v(y)) \frac{1}{|x-y|^{N+ 1}} (w(x)-w(y))\, dx\, dy\,.
\end{equation}
It is known that $\mathcal{B}(u,u)$ is an equivalent norm for the fractional Sobolev space $W^{1/2,2}(\ren)$. A number of basic properties of operator $\cal B$ are listed in \cite{CSV}, and we refer the reader to that reference when they are needed.

(ii) In our basic lemmas we use a variant where $f$ depends also on $x$ in a smooth way, see Subsection \ref{sec51}.1. The above inequalities are assumed to hold for such test functions.

\noindent  {\bf Note.} The solutions for which we derive the regularity result enjoy such properties and we call them {\sl weak energy solutions}. We will always assume that we are dealing with such class of weak solutions, i.e. they are bounded, integrable in space and obey the energy estimates.

%\newpage
%%%%%%%%%%%%%%%%%%%%%%%%%%%%%%%%%%%%%%%%%%%%%%%%%%%%%%%%%%%%%%%%%%%%%%%%%%%%%%%%%%%%%

\section{Strategy to prove H\"older regularity}\label{sec.holder}

Our concern in this paper is the local regularity of weak solutions.
The result we want to prove has been stated as Theorem \ref{mainthm}.
Let us review here the strategy and main tools. Since the equation is space- and time-invariant we may assume that $T_1<T_2 =0$, and then we may study the regularity around $x=0$ and $t=0$.
The basic ideas of the proof of regularity were proposed in \cite{CSV} and are as follows: On the one hand, we will prove some basic De Giorgi-type oscillation lemmas that say that the oscillation of the solution $u$ shrinks in a certain way when we restrict the domain with a certain scale. To be precise, we will rely on the iterated application of three basic lemmas, so-called {\sl oscillation reduction lemmas}. These technical results need only be proved for bounded nonnegative weak solutions defined in a strip $S_R=[-R,0]\times \RR^N$. We denote by $\Gamma_R$ the parabolic cylinder $[-R,0]\times B_R(0)$. By parabolic we mean at this point space-time subset. One of the lemmas controls the decrease  of the supremum  of the solution once we restrict the size of the parabolic neighborhood of \ $(0,0)$, from say $\Gamma_4$ into a smaller cylinder like $\Gamma_1$.
Another  one implies that under suitable assumptions the solution separates from zero in the same type of cylinders. A third one improves the first result so as to obtain a real alternative between going a bit down and a bit up. This is what is needed to make the iteration possible and efficient.

The possibility of efficient iteration depends on a second ingredient, the scaling property of the equation, that allows to renormalize the solution through the  transformation
\begin{equation}\label{scaling}
\widehat u(x,t)= A\, u(Bx, Ct)
\end{equation}
with $A,B>0$ free parameters, and  $C=AB$ (since $s=1/2$). Using this property, after every step of application of the oscillation reduction result, we renormalize the solution defined in $\Gamma_1$ into a scaled out solution defined in $\Gamma_4$ and we start a new application of the oscillation lemmas. In this way, we will show that the oscillation of the solution $u$ decays with a fixed geometric rate in a family of space-time cylinders shrinking also geometrically to a point.

\noindent $\bullet$ This program was  successfully implemented for $s<  1/2$ in paper \cite{CSV} but it turned out to be insufficient for $s\ge1/2$. A first difference in the situation  is that the reduction of oscillation of the basic lemmas (essentially, the first lemma) will also depend on the spatial $L^1$ norm of the solutions. Since this dependence gets worse  with the iterations because of the scalings, it has to be eliminated at some moment, and this happens through a geometrical transformation of the domain using a moving frame associated to some transport ideas that allow us to kill the integral term responsible for the unwanted dependence. This delicate idea has been successful for $s\in (1/2,1)$ as reflected in the final part of \cite{CSV}.

\noindent $\bullet$  However, the corrected method does not work for $s=1/2$ because the integral defining the change of coordinates, see formula \eqref{def.vel},  is in principle divergent (a divergence of logarithmic type) at both ends, zero at infinity, and the previous method corrects only the far field. We note that the problem arises only in the case where successive iterations push the solution always down so that we end up focusing at a point where $u$ vanishes (degenerate case). The intuitive reason for this difficulty to be anyway avoidable is that in such a situation the solution will be proved to be zero at some point with some H\"older exponent and this will make the problematic integral convergent at the origin. However, this is an argument a posteriori that has to be justified.

The technical way out of the difficulty that we propose here for $s=1/2$ consists of the consecutive applications of  scaling and a geometrical transformation, after separation of the effect of the velocity integral near the origin of coordinates (i.\,e., near the point under study). In the process a number of reductions of the domain are needed in a consecutive form, and they must be strictly controlled in order to arrive at the desired H\"older regularity estimate. The delicate analysis is explained in full detail in Sections \ref{s-large} and  \ref{s=12}, and it is the aim of the present paper.

 Before going into this  analysis we need to review the basic lemmas which support the proof, and we will devote the next sections to this task. We refer  to \cite{CSV} for more details on the derivation of the lemmas, which is done there for all $s\in (0,1)$. In the present shortened version we concentrate on the points of interest for the final sections.

\normalcolor

%%%%%%%%%%%%%%%%%%%%%%%%%%%%%%%%%%%%%%%%%%%%%%%%%%%%%%%%%%%%%%%%%%
\section{ The first oscillation reduction lemma}\label{sec51}

The first of the three basic lemmas deals with the question of ``lowering the maximum'' of a solution when we shrink the domain in a convenient way and the appropriate assumptions are met. It has been worked out in full detail in \cite{CSV} thinking of the case  $s< 1/2$. Special attention has to be paid to the difficulties of the case $s=1/2$ since the lack of integrability of the kernel at infinity means that we have to modify the argument. The statement is  formally almost the same as the lemma in \cite{CSV}, the difference lies in the additional dependence of the constants on $L^1$ information of the data. Recall the notations of Section \ref{sec.holder}. Our solution $u\ge 0$ is bounded above in the strip $S_4=\re^n\times [-4,0]$ in a precise way, and we must also assume that ${\cal M}_1=\sup\{\|u(\cdot,t)\|_{L^1(\ren)}: \ t\in (-4,0)\}$ is bounded. We write $\Gamma_4=B_4(0)\times [-4,0]$.

\begin{lem} \label{reg.1} Let $u \ge 0$ be an energy weak solution of \eqref{eqs} with $s=1/2$. Given $\mu\in (0,1/2)$ and  $\ve_0$  small enough,  there exists $\delta>0$ (depending possibly on $\mu,\ve_0, s, N$, and ${\cal M}_1$) such that, if we assume that \\ {\rm (i)}  the solution $u$ is bounded above in the strip $S_4=\re^n\times [-4,0]$ by
\begin{equation}\label{est.apriori1}
 \overline{\Psi}(x)= 1+(|x|^\ve-2)_+, \qquad 0<\ve <\ve_0\,,
\end{equation}
 and {\rm (ii)} $u$ is mostly below the level  $1/2$  in $\Gamma_4$ in the sense that
\begin{equation}
|\{u>1/2\}\cap \Gamma_4|\le \delta |\Gamma_4|\,,
\end{equation}
 then we have a better upper bound for $u$ inside a smaller cylinder:
 \begin{equation}
 \left. u\right|_{\Gamma_1}\le 1-\mu\,.
 \end{equation}
 \end{lem}

 Roughly speaking,  {\sl ``being mostly below 1/2 in space-time measure pulls down the supremum in a smaller nested cylinder''}, and this happens  in a quantitative form. Note that for this lemma $\delta$ can be chosen as a non-increasing function of $\mu$ with $\delta(1/2+)=0$.  Also, $\ve$ can be chosen as small as we want by sacrificing the gain in oscillation, i.\,e., the final H\"older exponent. We also remark that the size of the cylinders can be changed, though this affects the values of $\delta$ if the new sizes do not conform with the parabolic scaling. Finally, the levels $u=1/2$ and $u=1$ are taken
by convenience, any pair of levels $0<M_l<M_u$ will do, though in principle the value of $\delta$ will change.

\medskip

\noindent $\bullet$ Let us review the proof in order to recall the most important tools and formulas, and the final details that need attention. The basic idea in  the proof of the result is a particular kind of  ``localized energy inequalities'' that will be iterated
in the De Giorgi style  to obtain the reduction on the maximum in a smaller domain. Localization is obtained by using a suitable sequence of cutoff functions, which leads in the limit to the stated result.

In order to deduce the necessary energy inequalities we use integration by parts formulas and analysis of the  kernels. A main role is played by the bilinear form ${\cal B}(v,w)$ as defined above with kernel $K(x)=c|x|^{-(N+1)}$. Moreover, we put $L(x)=c_1|x|^{-N+1}$ so that $\Delta L =K$. We will repeatedly use the following equivalent form, based on \eqref{BB},
\begin{equation}
{\cal B}(u,v)=  \diint (u(x)-u(y))K (x-y) (v(y)-v(x))\,dxdy.
\end{equation}

\noindent {\bf \ref{sec51}.1.   An energy formula}.  We  consider a sequence of cutoffs $\vp(x)$ that have the form of downward perturbations of the level $u=1$  within a region containing the unit ball $B_1(0)$,  and they also have an  ``outer wing'' rising up above the 1-level for larger values of $|x|$ to be able to keep a global control of the dilations of $u$.  An explicit choice suitable for our purposes will be done below. We only need to know at this stage that the cutoff function $\vp$  is smooth, lies above $ 1/2$ everywhere, and also that $u\le \vp$ for all $|x|\ge 3$ for all times $-4\le t\le 0$.

We use the function $\eta=\log((u/\vp)\vee 1)=\log(g)$ as a  test function in the weak form of the equation, which is allowed in our definition of weak energy solution. Note that
$$
g:= \frac{u}{\vp}\vee 1 = 1+\frac{(u-\vp)_+}{\vp}=1+\frac{u_\vp^+}{\vp}\,,
$$
where $u_\vp^+=(u-\vp)^+$ according to the  notation introduced at the end of the Introduction. Note that $g\ge 1$ and $g>1$ iff $u>\vp$. According to our assumptions,  $u_\vp^+$ and $g-1$ have compact support in the ball of radius 3.  We will split $u$ as
$$
u=u_\vp^++\vp + (u-\vp)^-\,,
$$
where we write $(u-\vp)^-=(u-\vp)\wedge 0= u_\vp^-$. Notice that with this notation we have  $u_\vp^- \le 0.$  A detailed argument of paper \cite{CSV} shows that
 we have  the following basic identity for $T_1<T_2\le 0$:
\begin{equation}\label{energ.id}
\left\{\begin{array}{l}
\dint \vp (g\,\log g+1 -g)  \,|{_{T_2}} \,dx + \frac12\dint_{T_1}^{T_2}{\cal B}(u_\vp^+,u_\vp^+)\,dt+
\dfrac12 \dint_{T_1}^{T_2}{\cal B}(u_\vp^+,u_\vp^-)\,dt\\[10pt]
= \dint \vp (g\,\log g+1 -g) \,|{_{T_1}} \,dx -\dfrac12 \dint_{T_1}^{T_2}{\cal B}(u_\vp^+,\varphi)\,dt
+ \dint_{T_1}^{T_2}  {\cal Q}(u_\vp^+,u)\,dt\,,
\end{array}\right.
\end{equation}
where
$$
{\cal Q}(u_\vp^+,u)\,dt:=
\diint u_\vp^+(x)\frac{\nabla \vp(x)}{\vp(x)} \nabla_x L(x-y)[u(y)-u(x)]dxdy\,.
$$
 We will think of the LHS as the basic energy of this calculation, and the RHS as the terms still to be controlled.

\medskip

%%%%%%%%%%%%%%%%%%%%%%%%%%%%%%%%%%%%%%%%%%%%%%%%%%%%%%%%%%%%%%%%%%%%%%%%%%%%%%%%%%%%%%%%%%

\noindent {\bf \ref{sec51}.3. Cutoff functions, control of the RHS and final goal}

At this moment we make a convenient choice of the sequence of cutoffs in order to better tackle the RHS and rest of the proof of the lemma. Though only some simple bounds on the functions and their derivatives  are used, a practical choice used in \cite{CSV} is as follows:
 \begin{equation}
\vp_k(x)=\min\{1+(|x|^\ve-2)_+, \ \overline{\vp}_k(x)\}, \qquad \overline{\vp}_k(x) = \frac78 + \frac{|x|^2}{16}-\frac12 4^{-k}\,,
\end{equation}
for some small $\ve>0$ and $k=1,2,\dots$ Note that $\vp_k\ge \vp_{k-1}$. The following remark will be important: at points where $\vp_{k}<1$ we have
$$
\vp_k= \vp_{k-1}+\frac12 4^{-k}.
$$
We also have
$$
\inf \vp_k=\vp_k(0)> 1/2 \quad \mbox{  for } \ k\ge 1.
$$
Moreover, $\vp_\infty(x)\le 1$ precisely for $|x|\le \sqrt{2}$ and $\vp_1(x)< 1$ for $|x|< 2$. This means in particular that $\vp_k(x)=1+(|x|^\ve-2)_+$ for $|x|\ge 2$, $k\ge 1$.
Moreover, $\vp_\infty(x)=(|x|^2+14)/16 \le 15/16$ for $|x|\le 1$.

A more general version of the same construction takes
 \begin{equation}
\overline{\vp}_k= 1 -\frac1{2C} + \frac{|x|^2}{4C}-\frac12 C^{-k}
\end{equation}
with $C$ possibly larger than $4$. In that case $1-\vp_\infty(x)\ge 1/4C$ for $|x|\le 1$.

\medskip

\noindent  For the rest of this proof we write $u_k^+=(u-\vp_k)^+\ge0$, $u_k^-=(u-\vp_k)^-\le0$
with this choice of $\vp_k$. Notice that the support of $u_k^+$ is contained in the ball of radius 2 as a consequence of  assumption (\ref{est.apriori1}).

We are ready to tackle the RHS of  Identity \eqref{energ.id} with this choice of test functions. One part will be controlled by a small multiple of the present energy, i.\,e., we will absorb it into the LHS of \eqref{energ.id}. The rest will be bounded above by a large multiple of $|\{u_k^+>0\}|$ (we recall that the notation $|.|$ means the Lebesgue measure of the set).

\medskip

\noindent {\bf \ref{sec51}.4. Estimate of the remaining $\cal B$  term}

This part does not differ from \cite{CSV}. We start the process  with ${\cal B}_2={\cal B}(u_k^+,\vp_k).$ By inspecting the integral we easily get
$$
{\cal B}_2 \le \gamma {\cal B}(u_k^+,u_k^+)+\frac1{\gamma}{\cal B}^*(\vp_k,\vp_k)\,,
$$
for every $\gamma>0$, where ${\cal B}^*(\vp_k,\vp_k) $ indicates that the integral is performed only on the set where either $x$ or $y$ belong to  $\{u_k^+>0\}$. That is,
$$
{\cal B}^*=\diint [\chi_{\{u_k^+>0\}}(x) + \chi_{\{u_k^+>0\}}(y)  ]K(x-y)(\vp_k(x) -\vp_k(y))^2.
$$
For $\gamma $ small, then $\gamma {\cal B}(u_k^+,u_k^+)$ is absorbed into the LHS (into the energy). Now, using that
$$
|\varphi_k(x)-\varphi_k(y)| \le C \min (1, |x-y|),
$$
and the size of the kernel $K$, we arrive to the estimate
$$
{\cal B}^*
\le C |\{u_k^+>0\}|\le C 4^{2k}\dint_{\{u^+_k>0\}} (u^+_{k-1})^2dx
$$
 The last inequality follows by Chebyshev's inequality, since $u_{k-1}\ge 4^{-k}/2$ whenever $u^+_k>0$. The resulting expression is good for our later purposes.

\medskip

%%%%%%%%%%%%%%%%%%%%%%%%%%%%%%%%%%%%%%%%%%%%%%%%%%%%%%%%%%%%%%%%%%%%%%%%%%%%%%%%
\noindent {\bf \ref{sec51}.5. Analysis of the ${\cal Q}$ terms}

Let us finally examine the  last term in \eqref{energ.id}, a source of trouble for this paper. It also has a bilinear structure. Indeed,
$$
\begin{array}{l}
\diint u_k^+(x)\frac{\nabla \vp_k(x)}{\vp_k(x)} \nabla L (x-y)[u(x)-u(y)]\,dxdy := {\cal Q}(u_k^+,u) \\[14pt]
={\cal Q}(u_k^+,u_k^+) + {\cal Q}(u_k^+,\vp_k)+{\cal Q}(u_k^+, u_k^-)=
{\cal Q}_1+{\cal Q}_2+{\cal Q}_3\,.
\end{array}
$$
Note however that the ``kernel" that is involved is not symmetric due to the presence of terms with $\vp_k$. The study of the  contribution of each of the three terms is again split into the
close-range and far-field interactions, represented by the integrals where $(x-y)$ lies in a given ball around the origin, or alternatively in its complement. In that sense we note that  $\nabla L$ satisfies $|\nabla L|\le c|x-y|\,K(x,y)$. As in \cite{CSV}, the contribution of the first two terms in the RHS is either absorbed into the LHS or estimated
as smaller than
$$
C^{2k}\dint_{t_1}^{t_2}\!\!\dint_{\{u^+_k>0\}} (u_{k-1}^+)^2\,dx dt.
$$

\noindent $\bullet$ The last term in this analysis needs closer scrutiny since it is the source of the special difficulties for $s\ge 1/2$. We have
\begin{equation}\label{delicate}
{\cal Q}_3={\cal Q}(u_k^+, u_k^-)=-\diint u_k^+(x)\frac{\nabla \vp(x)}{\vp(x)} \nabla L (x-y)u_k^-(y)\,dxdy.
\end{equation}
 To estimate the integral on the set where $|x-y|\le \eta$ is small, we use that $|\nabla L (x-y)|\le C |x-y|\,K (x-y)$ and then $Q(u_k^+,u_k^-)$ is bounded by a small fraction of  ${\cal B}(u_k^+,u_k^-) $ (remember that we have proved that this term has the correct sign). We can therefore get this part absorbed  by the LHS of the energy identity.

\noindent $\bullet$ Finally, we need to consider the integral ${\cal Q}_3$ for $|x-y|>\eta $. This is the delicate case.Indeed, this is the only place where  the restriction $s<1/2$ was quite useful in our previous paper. We can solve the difficulty by making use of the known fact  that $u(x,t)$ is an  $L^1$ function in $x$, uniformly in $t$.  Since  $\nabla L(x-y)u_k^-(x) $ is bounded for large $|y|$, we get
$$
\begin{array}{l}
\left|\diint u_k^+(x)\frac{\nabla \vp(x)}{\vp(x)} \nabla L (x-y)\,u^-(y)\,dxdy\right|\\
\le C \|u(\cdot,t)\|_1 \left(\dint u_k^+(x)dx\right)\le C^{k}\|u(\cdot,t)\|_1 \left(\dint_{\{u^+_k>0\}} (u_{k-1}^+(x))^2dx\right).
\end{array}
$$
and we know that $\|u(\cdot,t)\|_1\le {\cal M}_1$. But we warn the reader that this idea is a partial solution and will run into difficulties later when we perform repeated iterations and rescaling.

\medskip

\noindent {\bf \ref{sec51}.6. Summary}. Putting all these estimates into \eqref{energ.id}, we  obtain  for $s=1/2$ (or any $s\in (0,1)$) and $t_1<t_2\leq 0$ \ the following energy inequality:
\begin{equation}\label{mod.ener.ineq}
\left   \{\begin{array}{l}
\dint \frac{(u_k^+ (t_2))^2}{\vp_k}dx + \frac12 \dint_{t_1}^{t_2}{\cal B}(u_k^+,u_k^+)\,dt\\[16pt]
\le
 2\!\dint \frac{(u_k^+(t_1))^2}{\vp_k }dx + (C^{2k} + C{\cal M}_1)\dint_{t_1}^{t_2}\!\!\dint_{\{u^+_k>0\}} (u_{k-1}^+)^2\,dx dt\,,
\end{array}\right.
\end{equation}
where $C$ is a universal constant that only depends  \normalcolor on  $s$ (here $s=1/2$) and the dimension $N$. We will have to pay attention to this dependence on $\|u(t)\|_1$ later, but in this first derivation of the lemma it allows us to continue and conclude. In the application to the iteration to follow next, the $t_i$ will be chosen in dependence of $k$.

\medskip

%%%%%%%%%%%%%%%%%%%%%%%%%%%%%%
\noindent {\bf \ref{sec51}.7. Iteration and end of proof of Lemma \ref{reg.1}}

This concluding argument follows the De Giorgi style and is done as in \cite{CSV}. Since this step is important and not long we recall it for convenience. We define the ``total energy function for the truncated solution'' $u_k^+$ as
\begin{equation}\label{def.aj}
{\cal A}_k=\sup_{ T_k\le t\le 0}\int  (u_k^+)^2(t)\,dx + \int_{T_k}^0 {\cal B}( u_k^+, u_k^+)\,dt,
\end{equation}
where $T_k=-2(1+2^{-k})$, $k=0,1,\cdots$. Notice that $\vp_k$ lies between $1/2$ and 1 at the points where $u_k^+$ is not zero. From (\ref{mod.ener.ineq}) with $k\ge 1$, and taking arbitrary values $t_2=t\ge T_k$ and $t_1=t' \in [T_{k-1},T_k]$ we have
\begin{equation} \label{bdd2}
{\cal A}_k\leq 4\inf_{t' \in [T_{k-1},T_k]} \int  (u_{k}^+)^2(t')\,dx+C^{2k}\dint_{t'}^{0}\!\!\dint_{\{u_k^+>0\}} (u_{k-1}^+)^2\,dx dt=I+II\,.
\end{equation}
Taking averages in $t'$ we arrive at the inequality
$$
\inf_{t' \in [T_{k-1},T_k]} \int  (u_{k}^+)^2(t')\,dx
\leq \frac 1{T_k-T_{k-1}}\int_{T_{k-1}}^{T_k} \int  (u_{k}^+)^2(t')\,dx dt'
$$
$$
\le 2^{k}  \int_{T_{k-1}}^{T_k} \int  (u_{k}^+)^2(t')\,dx dt'.
$$
Observing that $u_{k}^+(x)>0$ implies $u_{k-1}^+(x)>u_{k}^+(x)+4^{-k}/2$, we  realize that both, $I$ and $II$, have the same flavor, and that in fact we have the estimate
\begin{equation} \label{recur}
{\cal A}_k\leq C^{k}\dint_{T_{k-1}}^{0}\!\!\dint_{\{u_{k-1}^+>4^{-k}/2\}} (u_{k-1}^+)^2\,dx dt\,,
\end{equation}
for a possibly larger constant $C$.

The next step uses the following Sobolev embedding inequality
\begin{equation} \label{sob}
\left(\int u^p\,dx\right)^{2/p}\le C \|u\|^2_{H^{1/2}}
\end{equation}
for some $p>2$  depending on $r\in (0,1)$ and $N$. $C$ depends also on $r$ and $N$.
Actually, $p=2N/(N-1)>2$, when $N\ge 2$. Using this exponent and applying the inequality to $u_{k-1}^+$ we get
$$
\dint (u_{k-1}^+)^p\,dx \leq C \left[{\cal B}( u_{k-1}^+, u_{k-1}^+)\right]^{p/2}.
$$
Take $\theta=2/p$ and define $q=(1-\theta)2+\theta p$. Then
$$
\begin{array}{ll}
\dint_{\{u_{k-1}^+>4^{-k}/2\}} (u_{k-1}^+)^2\dx \le 4^{(k+1)(q-2)}\dint (u_{k-1}^+)^q\, dx \\
\le 4^{(k+1)(q-2)}
\left(\dint (u_{k-1}^+)^2\,dx\right)^{(1-\theta)}\left(\dint (u_{k-1}^+)^p\,dx\right)^{\theta}
\\
\le C  4^{k(q-2)} \left(\dint (u_{k-1}^+)^2\,dx\right)^{(1-\theta)}{\cal  B}(u_{k-1}^+,u_{k-1}^+)
\end{array}
$$
Integration in time $t$ along the interval $[T_{k-1}, \, 0]$ gives us from inequality \eqref{recur} and the previous estimate a recurrence relation  of the form
$$
\begin{array}{l}
{\cal A}_{k} \le   C^k\left( \sup_{ T_{k-1}\le t\le 0}\dint  (u_{k-1}^+)^2(t)\,dx\right)^{1-\theta} \cdot \dint_{T_{k-1}}^0 {\cal B}( u_{k-1}^+, u_{k-1}^+)\,dt \\
\le  C^k {\cal A}_{k-1}^{(1-\theta)}  {\cal A}_{k-1} =C^k {\cal A}_{k-1}^{1+\tau} ,
\end{array}
$$
with $\tau=1-\theta>0$ and a possibly larger constant $C$. It is well-known that this iterative sequence converges if ${\cal A}_1$ is small enough (depending on the constant $C$ appearing in the inequality. Applying this observation to our case we conclude that if we  take $\delta$ very small, then the iteration starts well so that the sequence ${\cal A}_{k}$ converges and then  ${\cal A}_\infty=0$, which means that $u\le \eta_\infty$ and this in turn implies that  $u\le 7/8$ for $|x|\le 1$. We thus get the result in the Lemma statement with $\mu=1/8$.

\medskip

 \noindent {\bf Remarks.} 1)  The obtained $\delta$ and $\mu$ depend on the bound  $\sup\{\|u(t)\|_1: t\in (-T,0)\}$. This dependence has to be eliminated later by a subtle new method.
 \normalcolor

 2) A simple modification of $\vp_\infty$ would give other values of $\mu\in (0,1/2)$, of course with a  different estimate of the maximum allowed value for $\delta$. The proof also shows that the time size $T=4$ can be replaced by any other number and the iteration will work with a different value for $\delta$ (and the same values for $\mu$ and $\ve$). \qed

%%%%%%%%%%%%%%%%%%%%%%%%%%%%%%%%%%%%%%%%%%%%%%%%%%%%%%%%%%%%%%%%%%
\section{The second basic lemma. Pulling up from zero}\label{sec53}

A similar oscillation reduction result applies from below. The proof is technically different since the equation is degenerate at $u=0$. The idea is that if $u$ is very often far from zero in $\Gamma_4$ then in a smaller, suitably nested cylinder  $u$ stays uniformly away from zero. The technical version explains how ``being  above $1/2$ most of the space-time, pulls  the solution up away from zero''.

\begin{lem} \label{reg.1b}  Under the  same assumptions set before Lemma \ref{reg.1}, given $\mu_0\in (0,1/2)$ there exists $\delta>0$ (depending possibly on $\mu_0,\ve_0, s$, and $N$) such that if  $u$ satisfies
\begin{equation}
|\{u \ge 1/2\}\cap \Gamma_4|\ge (1- \delta) |\Gamma_4| \,,
\end{equation}
 then \ $ \left. u\right|_{\Gamma_1}\ge \mu_0$.
\end{lem}

Again, $\delta$ is a non-increasing function of $\mu_0$. Let $\delta_0=\delta(1/4)$, that is, when $\mu_0=1/4$. It is important to remark that this result does not depend at all on the bound ${\cal M}_1=\sup\{\|u(\cdot,t)\|_{L^1(\ren)}: \ t\in (-4,0)\}$.

The proof of Lemma \ref{reg.1b} does not differ from the case $s<1/2$ done in \cite{CSV}.
A more elaborate version of this lemma will be needed in the second alternative of the iteration procedure.

%%%%%%%%%%%%%%%%%%%%%%%%%%%
%%%%%%%%%%%%%%%%%%%%%%%%%%%

\section{\bf The third oscillation reduction lemma}
\label{sec.3lemma}

 We still need to recall another ingredient before we attack the regularity issue by means of a suitable iteration. Indeed,  we have to improve Lemma \ref{reg.1} by showing that, in order to get a uniform reduction of the maximum in a smaller ball it is not necessary to ask that $u\le 1/2$ ``in most of'' $\Gamma_4$,  but only ``some of the time''. Most precisely, we must replace the sentence \lq\lq most of the space-time" of Lemma \ref{reg.1} by \lq\lq in some set of positive measure".

\begin{lem} \label{reg.2} {\rm (\lq\lq Some of the space-time below $1/2$, pulls down from 1'')} Assume as before that $0<s<1/2$ and  $u$ is trapped between $0$ and  $\overline{\Psi}$ in $S_4$. Besides, assume now that
\begin{equation}
|\{u<1/2\}\cap \Gamma_4|\ge \delta_0  |\Gamma_4|,
\end{equation}
with $\delta_0$ defined as above. Then $ \left. u\right|_{\Gamma_1}\le 1-\mu'$,  for some  $\mu'(\delta_0)$.
\end{lem}

Notice that this new lemma applies only in one direction, reducing the oscillation from above.
As in the classical porous media, we cannot expect this lemma to hold in the \lq\lq pulling-up" case, due to the property of finite propagation (existence of solutions with compact support), a consequence of the degeneration of the equation. Nevertheless, this one-sided improvement will be enough  to prove that the oscillation decays dyadically as explained shortly below.

The proof of this result is long and delicate but offers no difference with the one contained in Sections 9 and 10 of paper \cite{CSV}, using the so-called lemma on intermediate values. The precise version that is proved there is as follows.  Let us fix some notation. For $\lambda$ small enough w, we define for any $\ve>0$
$$
\psi_{\ve,\lambda}(x)=((|x|-1/\la^{4/s})^\ve-1)_+ \quad \mbox{if \ } \quad |x|\ge \lambda^{-4/s},
$$
and zero otherwise.

\begin{lem}\label{lem10.1} Given $\rho>0$ there exist $\ve>0$ and $\mu_1$ such that for any solution of the FPME in $ \ren\times(-4,0)$ satisfying
\begin{equation}
0\le u \le 1+\psi_{\ve,\lambda},
\end{equation}
and assuming that
\begin{equation}
|\{u<\vp_0\}\cap  (B_1\times(-4,-2))|>\rho\,,
\end{equation}
then we have
\begin{equation}
\sup_{ B_1\times (-1,0)} u \le 1-\mu_1.
\end{equation}
\end{lem}

 Note: In the next section we will take $\rho$ equal to $\delta_0$ as defined after the statement of  Lemma \ref{reg.1b}.

%%%%%%%%%%%%%%%%%%%%%%%%%%%%
%%%% OJO he leido hasta aqui
%%%%%%%%%%%%%%%%%%%%%%%%%%%%

\section{Iteration Procedure. Alternatives}\label{sec.iterproc}

The actual proof of the H\"older regularity result stated in Theorem \ref{mainthm} is based on
the application of the three basic lemmas, following the iterative process outlined in Section \ref{sec.holder}. We will review these steps only briefly since they have been explained in \cite{CSV}.  But the already announced difficulties motivate serious modifications that are the main contribution of the present paper. We will leave this delicate part for a later section.

The process works in an iterative way with two main alternatives. We want to  take any point $P_0=(x_0, t_0) \in \ren\times(0,\infty)$ and prove that $u$ is $C^\alpha$ around $P_0$
with an exponent that depends only on  $N$, and a H\"older constant that depends also on the $L^\infty$ norm of $u$ and a lower bound on $t_0$, $t_0\ge 4 \tau>0$.

Let us enter into some details: We have assume that our solution is bounded.  By scaling we may also assume  that $0\le u(x,t)\le 1$ in $\ren\times (0,T)$. Moreover, again by scaling we may assume that $T>t_0>5$. It will be then convenient to make a space-time translation and put $P_0=(0,0)$ assuming that the domain of definition of $u$ contains the strip $S_4=\ren\times [-4,0]$.

\noindent Consider now a positive constant $K < 1/4$ such that the growth of the outer wings is controlled as follows:
\begin{equation}\label{wing.growth}
 \frac1{1-(\mu_1/2)}\psi_{\lambda,\ve}(Kx)\le \psi_{\lambda,\ve}(x).
\end{equation}
 The coefficient $K$ depends only on $\lambda, \mu_1$ and $\ve>0$.
 The parameters are as in the last section. The iteration that we will perform offers two possibilities.

\medskip

 \noindent $\bullet$ {\sl Alternative 1. Regularity at a degenerate point.} Suppose that we can apply Lemma \ref{lem10.1} repeatedly  because the lowering of the oscillation may be assumed to happen always from above. We consider then the sequence of  functions  defined in the strip  $S_4=\ren\times (-4,0)$ by
  \begin{equation}\label{reg.scal}
 u_{j+1}(x,t)=\frac1{1-(\mu_1/4)}\, u_j(Kx, K_1 t), \quad K_1=\frac{K^{2-2s}}{1-(\mu_1/4)}\,.
 \end{equation}
 Note that this time the $u_j$'s are all of them solutions of the same equation.  According to the running assumption, and using \eqref{wing.growth}, we can apply Lemma \ref{lem10.1} at every step so that we have $u_{j}(x,t)\le 1-\mu_1$ in the cylinder $Q_1=B_1\times (-1,0)$ for every $j\ge 1$. In view of the scaling \eqref{reg.scal}, this would  imply H\"older regularity around the point  $(0,0)$, where the solution necessarily takes the degenerate value $u=0$ in a continuous way.
 This process was justified for $s<1/2$ in \cite{CSV}. But this process has a problem when we want to re-do the proof of the technical oscillation lemmas  for $s\ge 1/2$. In particular, the constants in Lemma \ref{reg.1} depend on the mass of the solution and this one grows unboundedly when performing successive rescalings, and this dependence propagates to the lemmas of Section \ref{sec.3lemma} that we are using. So the process deteriorates without control.

 \medskip

\noindent $\bullet$ {\sl Alternative  2. Regularity at points of positivity.} It can also happen that after some steps of the iteration the assumption on the measure of the set $\{u_j>1/2\}$ made in Lemma \ref{lem10.1}  fails. Then, we are in the situation where the oscillation is reduced from below thanks to Lemma \ref{reg.1b}, which pulls the solution uniformly up from zero in a smaller cylinder. Then the equation is no longer degenerate, because after that step we have
  $$
  0<\mu'\le u_j(x,t)\le 1,
  $$
in the cylinder $B_1\times (-1,0)$. Scaling the situation we will be in the conditions of the nondegenerate equation with  diffusivity $D(u)$ bounded above and below, so the case can be treated as quasi-linear. This case was carefully examined in the paper \cite{CSV}, where the proper modification of the proofs of the basic lemmas was discussed. In this way we obtain H\"older regularity at a point $P_0$ where $u(P_0)>0$ without modification on the arguments of the mentioned paper. We make some more detailed comments in Section \ref{sect.2alt}.

%%%%%%%%%%%%%%%%%%%%%%%%%%%%%%%%%%%%%%%%%%%%%%%%%%%%%%%%%%%%%%%%%%%%%%%%%%%%%%%%%%%%%%%%%%

\section{Correcting the iteration process for $s> 1/2$}\label{s-large}

As we have already indicated, Alternative 1 above  has a problem when we want to re-do the proof of the basic oscillation lemmas for $s\ge 1/2$. Indeed, we find a convergence problem  in the proof of Lemma \ref{reg.1}; the bulk of the proof contained in Section \ref{sec51} works without modification, and  an important difference was found only in the last estimate of Subsection \ref{sec51}.4, regarding integral ${\cal Q}_3$ in an outer region. The solution we have proposed in Subsection \ref{sec51}.4 was to make use of the extra fact  that $u(x,t)$ is also an  $L^1$ function in $x$, uniformly in $t$, and use this get the bound on the integral of the $y$ terms in ${\cal Q}_3$, with integrand $\nabla L(x-y)u_k^-(x) $, since it is bounded for large $|y|$.  This solves the problem of ending the proof of Lemma \ref{reg.1}, but then $\delta$ and $\mu$ depend on ${\cal M}_1=\sup_t \|u(t)\|_1$ as we have shown.

However, in order to obtain the $C^\alpha$ regularity result we have seen in the preceding Section \ref{sec.iterproc} that we need to iterate this (and the other oscillation  lemmas), we want to rescale and repeat, and then the difficulty re-appears, because we will keep stretching the variable $u$ and the $x$ axis, and therefore increasing the integral at every step, so the constants will be ruined in the iteration. We need a way
to control such  behaviour.

 It will be convenient to examine the whole part of ${\cal Q}(u^+,u)$ that contains the difficulty, i.e.,
\begin{equation}\label{q.bad.int}
   \diint u_k^+(x)\frac{\nabla \vp_k(x)}{\vp_k(x)} \nabla L (x-y) \,u(y)\,dxdy\,.
\end{equation}
 As we dilate and repeat in the iteration scheme, the term $ \nabla L(x-y) u(y)$
 also starts to build up as $y$ tends to infinity. On the other hand, the integrability in $y$ at infinity is lost if $1>s\ge 1/2$, since in that case $\nabla L$ decays  like
$$
 |\nabla L| \sim |y|^{-(N-2s+1)}
 $$
and this is not good enough. However, the good news is that for all $i,j$
    \begin{equation}\label{estim.diff}
 \partial_i \partial_j L \sim |y|^{-(N+2-2s)}\,,
\end{equation}
which is integrable as $|y|\to\infty$. Noting that $\nabla L$ is integrable for $y\sim 0$ if $s>1/2$, we conclude for such exponents that
$$
V(x,t): = (\nabla L (x-y)) \ast_y u(y,t)
$$
has a bounded H\"older seminorm. Therefore, it would be enough to control it at just a point, for instance at $x=0$ (for an interval of times). Let us see next how  this idea was implemented in \cite{CSV}.

%%%%%%%%%%%%%%%%%%%%%%%%%%%%%%%%%%%%%%%%%%%%%%%%%%%%%%%%%%%%%%%%%%%%%%
\subsection{ The transport approach for  $s>1/2$}

The technical way to make use of the last observation is to perform by a change of coordinates $x'=x-\gamma(t) $ that introduces a transport term to counter the difficult term we are dealing with, $ \int \nabla L(y)\,\Psi(y)\, u(y)\,dy. $ To be precise, we define
\begin{equation}\label{def.vel}
\gamma(t)=\dint_0^t \vec{v}(t)\,dt, \qquad
\vec{v}(t)=\dint \nabla L\, (y)\,u(y,t)\,dy\,,
\end{equation}
and we observe that $|\vec{v}(t)|$ depends also on $u$, and that  $|\vec{v}(t)|\le C<\infty $ since we are assuming that $u(y,t)$ is in $L^1_y$, uniformly in $t$. Indeed, the value of $\vec{v}(t)=\vec{v}(t; u)$ is only controlled by the integral of $u$ in space (what we call the mass of $u(t)$).
Next, we introduce change of variables
\begin{equation}\label{transf.v}
(x,t) \mapsto (x',t'):= (x-\gamma(t),t)\,,
\end{equation}
and we write the equation for $u$ with respect to the new variables, $u(x,t)=\tilde u(x',t')$. The RHS does not change since  we are performing a space translation for fixed time. However, the time derivative in the LHS transforms as follows:
$$
u_t(x,t) = \tilde u_{t'}(x -\gamma (t'),t')+(\nabla_{x'} u)(-\gamma'(t'))= \tilde u_{t'} - \vec{v}\cdot\nabla_{x'}   \tilde u.
$$
(We hope the reader will not have problems with the primes: in $\gamma'$ it means derivative, in $x', t'$ it means new space and new time). The last term in the formula is what we are aiming at.  The equation takes the convection-diffusion form
\begin{equation}\label{eq.diff.conv}
\tilde u_t-\vec{v}\cdot\nabla_{x'}   \tilde u = \nabla (\tilde  u\,\nabla {\cal L} \tilde u)\,.
\end{equation}
In the sequel we will write  $t$ for $t'$ and $u$ instead of $\tilde u$ without fear of confusion. The new space variable is still written $x'$. Next, we pass the term $\vec{v}\cdot\nabla u$ to the RHS and  multiply by $\log((u/\vp)\vee 1)$, as we did in Section \ref{sec51}, to obtain  the energy formula. We observe that in this case the RHS contains an extra term of the form
$$
I=- \diint  \nabla \log( (u/\vp)\vee 1) \,\vec{v}(t) u(x',t)\,dx'dt\,.
$$
This integral must be computed only in the region where $u>\vp$, and in that case $(u/\vp)\vee 1=u/\vp=1+(u_k^+/\varphi_k)$, so that
$$
I= -\dint dt\dint_{u>\vp_k}  \nabla u_k^+  \, \vec{v}(t) dx' +
\dint dt\dint_{u>\vp_k} u_k^+ \frac{\nabla \vp }{\vp}  \, \vec{v}(t) dx' =-I_1+I_2
$$
The first integral vanishes, and the second is precisely the troublesome  term:
$$
I_2=
\dint dt \dint_{u>\vp_k} dx \,\frac{u_k^+(x')}{\varphi_k(x')} \nabla \varphi_k(x')
 \dint \nabla L(y') u(y',t)\,dy'\,.
$$
 After this addition, the troublesome $\cal Q $ integral  (\ref{q.bad.int})   now amounts to
$$
\dint dt \dint_{u>\vp_k} dx \,\frac{u_k^+(x')}{\varphi_k(x')} \nabla \varphi_k(x')
 \dint (\nabla L(y')-\nabla L(y'-x')) \,u(y',t)\,dy'\,.
$$
 Using estimate (\ref{estim.diff}) this is convergent and can be estimated without having recourse to the $L^1$ norm of the solution.

Note.- The disappearance of the bad term  in the energy calculation in the new variables can be easily predicted if we write the equation for $\tilde u(x',t)$  in the more symmetrical form
\begin{equation}\label{trans.eq}
\tilde u_t+\nabla_{x'} \left(\tilde u  \left(\int \{\nabla L(y-x')-\nabla L(y)\}(\tilde u(y)-\tilde u(x'))\,dy\right)\right)=0\,,
\end{equation}
to be interpreted in the same weak form, or weak energy form, as  we have used for $u(x,t)$.

In any case, this allows  to prove Lemma \ref{reg.1} also for $s\in (1/2,1)$ if we work in the new coordinates, and the constants involved in the result do not depend on the $L^1$ norm of the solution. The price to pay is that the slope of the distorted space variables does depend on the $u$-integral. So, in the first step of the iteration process we have shown how to transfer the difficulty from a numerical term to a geometrical distortion.

In order to sum up the result, let us introduce the bound $M=1\vee \sup_{t>0} \vec{v}(t)$,
that depends only on $u$ via the norm $\sup_t\|u(\cdot,t)\|_1$.

\begin{lem} \label{reg.1Mod} Let $1/2< s<1$ and let $u$ be a solution of the FPME under the assumptions of Lemma {\rm \ref{reg.1}}. Let us perform the above change of variables so that $\tilde u(x',t')$ is defined in
smaller cylinder $Q_L$ where $L=4/(M+1)$. Then the result of Lemma {\rm \ref{reg.1}} is true for $\tilde u$, with conclusion holding in a smaller cylinder $Q_{1/M}$;  $\delta$ may depend also on $M$.
\end{lem}

Thanks to \eqref{trans.eq},  it is then immediate to see that the pull-up Lemma \ref{reg.2} are also true if stated in the form that we have used for Lemma \ref{reg.1Mod}. A bit more of attention to the details will show that the stronger reduction Lemma \ref{lem10.1} also holds, since the iterations do not change the scaling in space and time.

%%%%%%%%%%%%%%%%%%%%%%%%%%%%%%%%%%%%%%%%%%%%%%%%%%%%%%%%%%%%%%%%%%%%%%
\subsection{Analysis of the transport term in the iterations}

When we try to perform again the iteration procedure of Section \ref{sec.iterproc}, one of the alternatives is repeated scaling around a degenerate point. In that case the iterations take the form
\begin{equation}
 u_{j+1}(x,t)=\frac1{1-(\lambda_*/4)}\, u_j(Kx, K_1 t), \quad K_1=\frac{K^{2-2s}}{1-(\lambda_*/4)}\,.
 \end{equation}
 that we may sum up as
$$
u_1(x_1,t_1)=A\,u(x,t) , \quad x_1= B x, \quad t_1=C t.
$$
where $A>1$, $B<1$ and $C=B^{2-2s}A$, so that the same equation will be satisfied after the
change of scale. We propose here to do the same iteration for the solution $\tilde u$ in terms of the variables $x'$ and $t$. The equation will then take the modified form \eqref{trans.eq}, that will be satisfied again by the iterates, just as it is written.  It is true that the velocity $\vec{v}(t)$ will change from iteration
to iteration according to the rule
\begin{equation}\label{vel.change}
\vec{v_1}(t)= \frac{C}{B}\,\vec{v}(Ct)=B^{1-2s} \vec{v}(Ct),
\end{equation}
which follows both from the geometric transformation, and from the definition of $\vec{v}$ in \eqref{def.vel}.
Therefore, after the first geometrical transformation such repeated iterations  conserve the same correspondence for all subsequent steps. In other words, the geometrical transformation done in the first step will hold for all remaining steps: if the set of coordinates at that moment is  $(x_n,t_n)$, we obtain a set of newly distorted  coordinates $(x'_{n},t_n)$ by the formula
$$
x'_{n}(t)=x_{n}(t)-\gamma_n(t), \qquad \gamma_v'(t)=\vec{v}_{n}(t)
$$
which is just a scaled version of the original transformation for $n=0$. Summing up, since the contractions
in the upper bound for $u$ happen with a constant rate $1-\mu$ in cylinders that shrink in space and time also  with a fixed rate, we conclude in a standard way $C^\alpha$-regularity with respect to the transformed variables. But since the coordinate  transformation is done only once and is Lipschitz continuous, this means the same type of H\"older regularity for $u$ with respect to the original coordinates $(x,t)$. Of course, the Lipschitz constant of the transformation depends on ${\cal M}_1=\sup_t \|u(t)\|_{L^1_x}$.

%%%%%%%%%%%%%%%%%%%%%%%%%%%%%%%%%%%%%%%%%%%%%%%%%%%%%%%%%%%%%%%%%%%%%%%%%%%%%%%%%%%%%%%%%%
%%%%%%%%%%%%%%%%%%%%%%%%%%%%%%%%%%%%%%%%%%%%%%%%%%%%%%%%%%%%%%%%%%%%%%%%%%%%%%%%%%%%%%%%%%

\section{The transport approach for $s=1/2$}\label{s=12}

The problem with the $\cal Q$ estimate in Lemma \ref{reg.1} was solved in Section \ref{s-large}
by means of the change of variables described as transport approach. It allowed us to swallow the conflicting term by means of  a controlled geometrical distortion on the assumption that $s\in (1/2,1)$. Unfortunately, this does not work  for $s=1/2$ because the integral defining the relative velocity of the new coordinates in \eqref{def.vel} is now given by the expression
\begin{equation}\label{def.vel.12}
\vec{v}(t)=\dint \nabla L\, (y)\,u(y,t)\,dy= C\,\dint \frac {u(y,t)\,\sign(y)}{|y|^{N}}\,dy\,,
\end{equation}
which is not necessarily bounded. It would be at $t=0$ if we already knew that $u(0,0)=0$ and moreover that there is a modulus of continuity of $u(0,x)$ in $x$. But this is precisely what we want to prove. Therefore, we have to work in an indirect way, by using partial transport corrections in an iterative way. This happens as follows:

%%%%%%%%%%%%%%%%%%%%%%%%%%%%%%%%%%%%%%%%%%%%%%%%%%%%%%%%%%%%%%%%%%%%%%%%%%%%%%%%%%%

\subsection{Transport after isolating the origin}
We go back to the beginning of the transport method, as introduced in the previous section, but now we make a  partition
$$
 \nabla L= (\nabla L) \Psi+(\nabla L )({1-\Psi})\,,
$$
with  a smooth cutoff function $0\le \Psi\le 1$ such that $\Psi\equiv 0$ in $B_\ve$ and $\Psi\equiv 1$   out of    $B_{2\ve}$, i.\,e., it eliminates a neighborhood of the origin. The term corresponding to  $\nabla L (1-\Psi)$  is controlled by a small multiple of the ``good term"
$$
\diint u^+_\vp(x)K(x-y)(u-\vp)^-(y)\,dxdy
$$
since  $\nabla L (1-\Psi)=0$ outside $B_{2\ve}$ and  $\nabla L\,(1-\Psi)\le |x-y|K(x-y)\le \ve K(x-y)$.
We now consider the term $\nabla L\,\Psi$. Again, we observe that
\begin{equation}\label{ijestimate}
 \partial_i \partial_j  L \Psi \sim \frac1{(1+|y|)^{N+1}}
\end{equation}
(valid for all second derivatives), hence
$$
V(x,t): = (\nabla L\Psi(x,y,t)) \ast_y u(y,t)
$$
has a bounded Lipschitz seminorm. It will be enough to control it at just a point, for instance at $x=0$ (for an interval of times), since the spatial increments are already under control.

We can now perform the coordinate change
$x'=x-\gamma(t) $ that introduces a transport term to eliminate the term $ \int \nabla L(y)\,\Psi(y)\, u(y)\,dy. $ To be precise, we define
\begin{equation}\label{def.vel2}
\gamma(t)=\dint_0^t \vec{v}(t)\,dt, \qquad
\vec{v}(t)=\dint \nabla L\,{\Psi} (y)u(y,t)\,dy\,,
\end{equation}
and we observe that $|\vec{v}(t)|\le C<\infty$, since $u(y,t)$ is in $L^1_y$ uniformly in $t$. Indeed, the value of $\vec{v}(t)=\vec{v}(t; u,\Psi)$ is only controlled by the integral of $u$ in space (i.\,e., the mass of $u(t)$).
When we make the change of variables
$$
(x',t)= (x-\gamma(t),t)\,,
$$
and we write the equation with respect to the new variables, this  does not change the RHS. However, on the LHS we get
$$
u_t(x,t) = \tilde u_t(x-\gamma (t),t)-(\nabla \tilde u_{x'} u)\,\gamma'(t)= \tilde u_t - \nabla_{x'} \tilde u\cdot\vec{v}\,.
$$
The equation becomes
\begin{equation}\label{trans.eq2}
\tilde u_t+\nabla_{x'} \left(\tilde u  \left(\int \{\nabla L(y-x')-\nabla L(y)\Psi(y)\}\tilde u(y)\,dy\right)\right)=0\,,
\end{equation}
so that the kernel $L(x-y)$ has been replaced by $\tilde L(x,y)=L(x-y)+L(y)\Psi(y)$. Let us see how this affects the energy estimates. If we pass $\vec{v}\cdot\nabla u$ to the RHS and  multiply  by $\log((u/\vp)\vee 1)$, as we did above to obtain in the energy formula, this one will contain an extra term of the form
$$
-\dint dt \dint dx \nabla \log( (u/\vp)\vee 1) \,\vec{v}(t) u(x)=I.
$$
The integral must be performed only when $u>\vp$ an in that case $(u/\vp)\vee 1=u/\vp=1+(u_k^+/\varphi_k)$ so that
$$
I= -\dint dt\dint_{u>\vp_k}  \nabla u_k^+  \, \vec{v}(t) dx+
\dint dt\dint_{u>\vp_k} u_k^+ \frac{\nabla \vp }{\vp}  \, \vec{v}(t) dx =-I_1+I_2
$$
The first integral vanishes and the second is the desired term:
$$
I_2=
\dint dt \dint_{u>\vp_k} dx (u_k^+(x)/\varphi_k(x)) \nabla \varphi_k(x)
 \dint \nabla L(y)\,{\Psi} (y)u(y,t)\,dy\,.
$$
 The troublesome integral in $\cal Q $  now becomes
$$
{\cal Q}^*= \dint dt \dint_{u>\vp_k} dx \,\frac{u_k^+(x')}{\varphi_k(x')} \nabla \varphi_k(x')
 \dint [\nabla L(y')\Psi(y')-\nabla L(y'-x')] \,u(y',t)\,dy'\,,
$$
which is estimated in view of the decay rates (\ref{ijestimate}) of the second derivatives of $L$,  and the fact that $x'$ is bounded.  This correction simplifies the energy formula and allows to prove the First Lemma in the new coordinates with constants that do not depend on the $L^1$ norm of the solution. The price to pay is that the slope of the distorted space variables does depend on the sup of the $L^1$ integral, as we have seen in the previous section. So in the first step of the iteration it seems that this argument does not any definitive improvement since all amounts to an estimate with worse constants.

%%%%%%%%%%%%%%%%%%%%%%%%%%%%%%%%%%%%%%%%%%%%%%%%%%%%%%%%%%%%%%%%%%%%%%55
\subsection{The iteration. First step}
 It is convenient at this stage to take some notational steps  since we are going to produce a sequence of solutions $u_j$ in a nested sequence of domains $\widehat Q_j\subset \widehat S_j$, but the estimates are going to be obtained after successive changes of coordinates, and we need to carefully label them and their domains. We will use dots for new coordinates, while the derivative in time of a function  $\gamma(t)$  is represented by $\dot{\gamma}(t)$.

We start with the standard cylinder $\widehat Q_1=Q_4(0,0)$ included in the strip $\widehat S_1:=S_4=\ren\times (-4,0)$ and we label this original solution as $u=u_1$. We assume that $u_1\le 1$ in $\widehat Q_1$.

We introduce the first correcting speed $\vec{v}_1(t)$, which is bounded by $C_1$ (that depends on the sup of the $\|u(t)\|_1$). We perform the  change of coordinates defined by $x'=x- \gamma_1(t)$  with $\dot{\gamma}_1(t)=\vec{v}(t)$, given by \eqref{def.vel2}. Time in not changed in this step. The resulting space-time transformation is denoted by ${\cal T}_1=\widehat Q_1\mapsto \widehat S_1$. Due to the convective effect, it will happen that ${\cal T}_1(\widehat Q_1)$ is not contained in $\widehat Q_1$, hence we restrict the domain to a reduced cylinder $\widehat Q_{1,r}$ (subindex $r$ for reduced) in such a way that ${\cal T}_1(\widehat Q_{1,r})\subset \widehat S_1$. Indeed, the target  domain can be chosen to contain  $\widehat Q_1'=Q_{4/C'_1}$ with $C'_1>1$ depending only on $C_1$. We denote by $u_1'(x',t)=u_1(x,t)$ the transformation of function $u_1$ by ${\cal T}_1$ (primes indicate here new coordinates).

We can now apply  the first and third lemmas to the modified equation satisfied by $u_1'$ to get a reduction of the upper bound of $u_1'$ in a quarter domain $\widehat Q_1'/4$. We get $u_1'\le 1-\mu$ in $\widehat Q_1'/4$. Note that this means that $u_1(x,t)\le 1-\mu$ in $\widehat Q_{1,s}= {\cal T}_1^{-1}(\widehat Q_1'/4)$, which has the shape of a cylinder with base  $B_{1/C_1}$; however, it is not a vertical cylinder but one with curved director line slanted in the direction of the line $x=\gamma_1(t)$ with $\gamma_1'(t)= \vec{v}_1(t)$. See Figure 1.

%Los dos dibujos en una linea

\begin{figure}  \label{fig.1}
    \includegraphics[height=6cm, width=6.1cm]{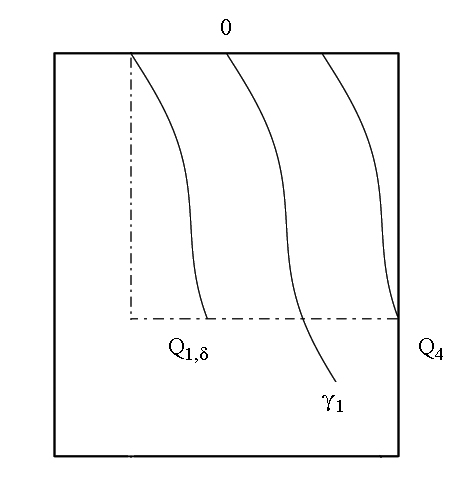}
\qquad\qquad
    \includegraphics[height=6cm, width=6.1cm]{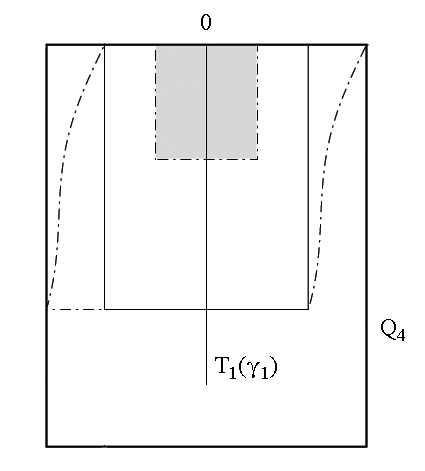}
 \caption{{\textit{Sketch for the First Step of the Iteration}}
 \newline
\qquad \qquad \qquad \textit{Left: original coordinates. \qquad Right: coordinates after the transformation}
 }
 \end{figure}

%%%%%%%%%%%%%%%%%%%%%%%%%%%%%%%%%%%%%%%%%%%%%%%%%%%%%%%%%%%%%%%%%
\subsection{Extension between iteration steps}

We now make a scaling of the variables (or extension, ${\cal E}_1$) that passes from $u_1(x,t)=u(x,t)$ \  to \ $u_2(x_2,t_2)$ of the form
\begin{equation}\label{scal.12}
u_1(x, t)=Au_2(x_2,t_2), \quad x_2=Bx, \quad t_2=T t\,,
\end{equation}
with parameters $A,B$ and $T>0$. We take $A=1-\mu$ to compensate for the shrinking gain of the previous step, while $B>1$  will be chosen to recover more or less the  cylinder $Q_4$; we finally need to put $T=BA$ \ so that the equation  satisfied by $u_2$ w.r.t $(x_2,t_2)$ will be the same FPME \eqref{eq1} satisfied by $u(x,t)$. Then, we can easily check that the old speed and the new speed (defined as in \eqref{def.vel2}) are related by
\begin{equation}
\vec{v}(t; u_1,\Psi)= A\,\vec{v}(t_2;u_2,\Psi_1)\,, \qquad \Psi(y)=\Psi_1(By)\,,
\end{equation}
where $t_2=Tt$. This relation also follows from transforming equation \eqref{eq.diff.conv}. Since $A=1-\mu<1$, the velocity increases slightly in norm in the extension transformation. This factor $A<1$ will play an important role later.

The scaling can also be applied after transformation ${\cal T}_1$, thus using  the variables $x_2'=Bx'=x_2-B\,\gamma_1(t)$ and
$$
u_1'(x',t)=Au_2'(x_2',t_2)
$$

When we perform the same type of the energy and transport estimates after this scaling we will get a remaining term (to be controlled) of the form
$$
{\cal I}_2(t_2)=\dint \nabla L (z) {\Psi}(z)\, u_2(t_2,z)\,dz\,,
$$
and now $u_2$ has a larger integral, $\|u_2(t_2)\|_1=(B/A)\|u_1(t)\|_1$. This allows to calculate the velocity $\vec{v}_2(t_2)$ of the new coordinate transformation. Note that the corrected $L$ of the first step transforms into
$$
\tilde L(x,y)=L(x-y)+L(y)\Psi(y)=B^{-(N-1)}(L(x_2-y_2)+L(y_2)\Psi(y_2/B))
$$
The presence of the last term means that we that we have already corrected from the $\cal Q$ energy integral an amount equivalent to the relative velocity
$$
\vec{v}_1(t; u,\Psi)=\int \nabla L (y) {\Psi}(y) u(t,y)\,dy=A \int \nabla L (y_2) {\Psi}\,(y_2/B) \, u_2(t_2,y_2)\,dy_2\,.
$$
In conclusion, in the $x_2'$ coordinates we need only be concerned about compensating for the remainder
\begin{equation}
\vec{v}_{12}(t_2)= \dint \nabla L(y_2)\,(\Psi(y_2)-\Psi(y_2/B))\,u_2(t_2, y_2)\,dy_2\,,
\end{equation}
with cutoff function $\Psi_{12}(z)=\Psi(z)-\Psi(z/B)$, which reflects the correction to pass from one iteration to the next. Since $\Psi_{12}$ is supported in the annulus $1\le |z|\le 2B$,  \sl we can estimate this new  correcting speed in a uniform way without using any information on the integral of $u$ in the whole space, only the local $L^\infty$ estimates produced in the iteration. \rm This is the crucial observation for the success of our plan.

%%%%%%%%%%%%%%%%%%%%%%%%%%%%%%%%%%%%%%%%%%%%%%%%%%%%%%%%%%%%%%%%%%%%%%
\begin{figure}  \label{fig.2}
    \includegraphics[height=6.0cm, width=6.2cm]{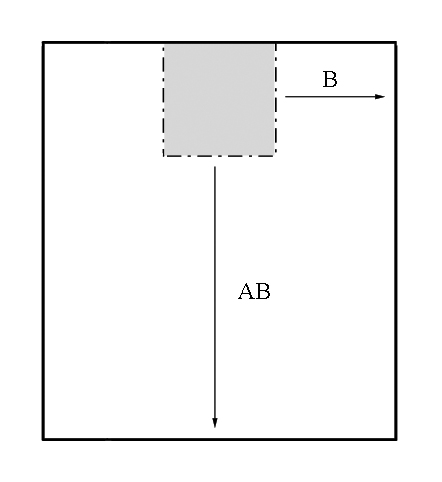}
\qquad
    \includegraphics[height=7cm, width=7.3cm]{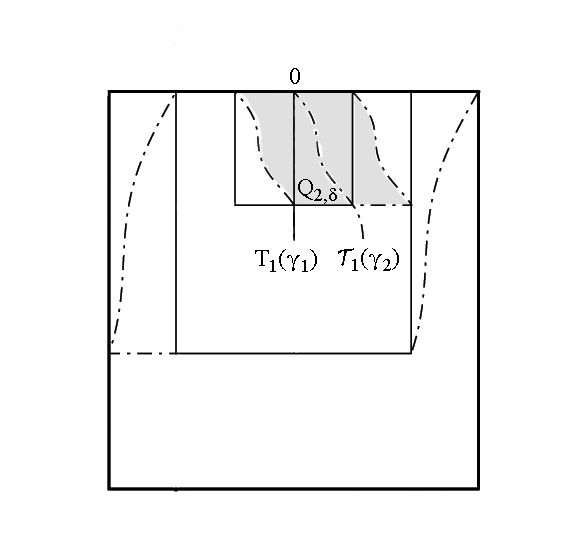}
 \caption{\noindent\textit{Left: sizes of the extension.} \textit{Right: Starting the second step of the iteration}}
\end{figure}

\subsection{Second step of the iteration: new reduction }

 We proceed with the new correction.  We pass to a new moving frame, a further correction to the transformation we have already done to produce $u_1'$ from $u_1$, so that we can apply the reduction lemmas. We thus put
$$
u_2'(x_2',t_2)=u_2(x_2,t_2)\,,
$$
and $x_2'$ is obtained from $x_2$ after the change of variables  that takes into account the
velocity $\vec{v}_2$. Let us call this correspondence $T_2$. Now
$$
\frac{dx_2'}{dt_2}-\frac1{A}\frac{dx'}{dt}=v(t_2,u_2,\Phi)-v_2(t_2,u_2,\Phi(y/B))=v_{12}(t_2).
$$
The correction in the slant of the image cylinder when referred to the standard scaling of $\widehat Q_1'/4$ is thus the relative speed $\vec{v}_{12}$ that is uniformly bounded, even the original speed $\vec{v}_1$ has been amplified, but we need not worry about it. We make the corrected change of variables ${\cal T}_2$ with suitable cylinders that shrink from $Q_4$ in a controlled way. The starting estimate is $u_2'\le 1$ in that cylinder.

We can now apply  the first and third lemmas to get a reduction of the upper bound of $u_2'$ in a quarter domain $\widehat Q_2'/4$. We get $u_2'\le 1-\mu$ in $\widehat Q_2'/4$. Note that this means that $u_2(x,t)\le 1-\mu$ in $\widehat Q_{2,s}= T_2^{-2}(\widehat Q_2'/4)$, which has the shape of a cylinder with base a ball, but not a vertical cylinder since it has a curved director line slanted in the direction of the line $x=\gamma_2(t)$ with $\gamma_2'(t)= \vec{v}_2(t)$.

In the non-scaled domains, i.\,e., after undoing  scaling \eqref{scal.12k} and transformation ${\cal T}_1$, the estimate piles up two contractions $(1-\mu)^2$ in a smaller slanted domain, whose inclination is $v_1$ plus $Av_{12}$. Observe the factor $A<1$ in the new term, derived by taking advantage of the size reduction of the first step; this is what will help us in summing the iterations.

%%%%%%%%%%%%%%%%%%%%%%%%%%%%%%
\subsection{Further steps and conclusion of the first alternative}
New correcting speeds can be calculated for the subsequent iterations and all of them will be bounded by the same constant when calculated in the domain obtained after scaling and passing to the corresponding moving frame. The scaling has the form
\begin{equation}\label{scal.12k}
u_k(x_k, t_k)=Au_{k+1}(x_{k+1},t_{k+1}), \quad x_{k+1}=Bx_{k}, \quad t_{k+1}=T t_{k}
\end{equation}
with $A=1-\mu$, $B>1$  and  $T=BA$ as before.  We will perform in every step of the iteration the corresponding change of variables immediately after the scaling. If the set of coordinates at that moment is  $(x_k,t_k)$, we obtain a set of newly distorted  coordinates $(x'_{k},t_k)$ by the formula
$$
x'_{k}(t)=x_{k}(t)-\gamma_k(t_k), \qquad d\gamma_k/dt_k=\vec{v_k}
$$
but the change with respect to applying scaling to the previous moving frame will be a correction speed $\vec{v}_{k-1,k}$ that is uniformly bounded as we have shown in the second step. We can then reduce the domain by a certain factor $C_1$ to make it fit into the distortion of the standard $Q_4$, apply the reduction lemmas and conclude that in a much smaller domain (but always the same proportion) we get $u_k(x_k,t_k)\le 1-\mu$.

Undoing the $k$ scalings and $k$ relative transformations we arrive at a small slanted cylinder in the original variables $\widehat Q_k$, with sizes $L_k\sim B^{-k}$ in space and $T_k\sim (AB)^{-k}$ in time. In such sequence of cylinders the solution $u(x,t)$ satisfies $u\le (1-\mu)^k$ and
the slant has accumulated a speed that is bounded by
$$
C_k=C_1+\sum_1^k CA^k
$$
and only $C_1$ depends on an $L^1$ norm. We conclude that $C_k$ is uniformly bounded as $k\to\infty$.

In this way we conclude that at a degenerate point where $u(x_0,t_0)$ the solution
has the $C^\alpha$ estimates in a backward parabolic neighborhood with an $\alpha$ that depends on the constants that we have been carrying around, in the end functions of $N$. Of course, the constant in the H\"older seminorm depends on the $L^\infty$ and $L^1$ norms of the solution in the strip $S$
where we do the calculations, in the end on the $L^1$ norm of the initial data, but this can be derived from the scaling group of the equation. This is what gives the $C^\alpha$ regularity.

\medskip

\noindent {\bf Comment.} We have performed a very detailed, step by step construction of the process, but the crux of the new argument can be expressed in simpler words. The convergence of the corrections relies on the possibility of finding a sequence of cylinders that shrink more in one direction that in the other (less in time than in space) in the original coordinates. This has to be combined with a sequence of geometrical distortions based on iterated cutoffs. Thus, we are keeping the same ratio of dilations in all the iterations, but the H\"older continuity could be obtained even if we relax  that requirement on the condition of keeping a geometrical series behaviour and a longer time than space scale.

%%%%%%%%%%%%%%%%%%%%%%%%%%%%%%%%%%%%%%%%%%%%%%%%%%%%%%%%%%%%%%%%%%%%%%%%
\section{The second alternative}\label{sect.2alt}

 Another alternative in the iteration procedure of Section \ref{sec.iterproc} happens when in one of steps we reduce the oscillation by below by using Lemma \ref{reg.1b}, so that we pull the solution away from zero. Therefore, we get a situation where the solution is bounded between two positive constants $0<M_l<M_u$ in the reference cylinder.  The technical details of Alternative 1 work the same way in this
case and we only have to choose the sense of the oscillation reduction,
applying it from above or below, depending on the solution being above or
below the middle value of the strip most of the time. Note that in Lemma
\ref{reg.1} the solution never sees the degenerate part of the equation, due to
the nature of the cut off.

Indeed, this case is easier since from the moment we apply the lemma where the solution is pulled up we have uniform ellipticity in the iterations, we will be even converging along the iterations to an equation with constant diffusivity coefficient.
This variant of the analysis has also been discussed  in detail in the paper \cite{CSV}, and the precise modification has been performed in Sections 7 of that paper, hence we may  leave the  details to the reader. This ends the proof of Theorem \ref{mainthm} for $s=1/2$, after recalling as technical tool the modification of the energy calculation performed in Section 7 of \cite{CSV}.

%%%%%%%%%%%%%%%%%%%%%%%%%%%%%%%%%%%%%%%%%%%%%%%%%%%%%%%%%%%%%%%%%%%%%%%%%%%%%%%%%%%%

\section{The one-dimensional case}\label{sect.1d}

We have assumed throughout the paper $N\ge 2$ since the case $N=1$ has some peculiarities worth commenting. As we have already indicated in our first paper \cite{CaVa09}, the functional treatment of the equation and its kernels is very similar for $s<1/2$ but for $s=1/2$ we find that $L(x-y)=c\log(x-y)$ which is unbounded at infinity. The way out of the difficulty is to avoid the consideration of the pressure and work always with the derivatives $\nabla_x L(x-y)\,u(y)$ that actually appear in the equation, and involve no growth at infinity. This implies some work in revising all the proofs above, but the main result will still hold.

On the other hand, the 1-dimensional theory has alternative existence proofs and very good extra mathematical properties. Thus,  Biler, Karch and  Monneau study this equation in \cite{BKM} as a model for dislocation and find that the integrated version admits viscosity solutions that are unique and admit comparison. On the other hand, Carrillo, Ferreira and Precioso \cite{CFP2012} apply transportation methods and show that the solution can be obtained as a gradient flow in the space ${\cal P}_2$ of probability measures with bounded second moment, which implies that the maps $u_0\mapsto u(t)$ form a contraction semigroup. Such properties are not proved in the multidimensional case $N\ge 2$.

%%%%%%%%%%%%%%%%%%%%%%%%%%%%%%%%%%%%%%%%%%%%%%%%%%%%%%%%%%%%%%%%%%%%%%%%%%%%%%%
\section{Extension of the existence theory}\label{sec.ext}

After these results, we can extend the existence theory to all nonnegative and integrable initial data. This was already consider in \cite{CSV} for $s\ne 1/2$.

\begin{thm}\label{thm.existence} For every $u_0\in L^1(\ren)$, $u_0\ge 0 $, there exists a continuous weak solution of the FPME \eqref{eq1} in the following sense: there exists a function $u(x,t)$, continuous and nonnegative in $Q=\RR^n\times (0,T)$ such that \
$$
u\in L^\infty(0,\infty:L^1(\ren)\cap L^\infty(\ren\times (\tau,\infty) \ \mbox{ for all \ } \tau>0\,,
$$
$$
{\cal K}(u)\in L^1(0,  T: W^{1,1}_{loc}(\RR^n)), \qquad u\,\nabla{\cal K}(u)\in L^1(Q_T)
$$
and the identity
\begin{equation}\label{identity}
\iint u\,(\eta_t-\nabla {\cal K}(u)\cdot\nabla\eta)\,dxdt+ \int
u_0(x)\,\eta(x,0)\,dx=0
\end{equation}
holds for all  continuously differentiable test functions $\eta$ in $Q_T$  that are compactly supported in the space variable and vanish near $t=T$.
\end{thm}

The proof does not depart from the one performed in \cite{CSV}. Taking this proof into account we may eliminate the boundedness condition from the assumptions on the solutions in our  main theorem by just asking the the initial data are integrable. On the other hand, once we have a theory for integrable data we can extend it to nonnegative Radon measures as initial data as in done in \cite{SerVaz}.

\

%%%%%%%%%%%%%%%%%%%%%%%%%%%%%%%%%%%%%%%%%%%%%%%%%%%%%%%%%%%%%%%%%%%

\noindent {\large\bf Acknowledgments}.   L. Caffarelli has been funded by NSF Grant DMS-0654267  (Analytical and Geometrical Problems in Non Linear Partial Differential Equations),  and J. L. V\'azquez by Spanish Grant MTM2011-24696. Part of the work was done while both authors were visitors at the Isaac Newton Institute, Cambridge, during the Free Boundary programme 2014.  We thank the support and hospitality. %Graphs were prepared by M. Bonforte.

\vskip 1cm

%\newpage

%%%%%%%%%%%%%%%%%%%%%%%%%%%%%%%%%%%%%%%%%%%%%%%%%%%%%%%%%%%%%%%%%%%%%
\bibliographystyle{amsplain} 

\begin{thebibliography}{10}


\bibitem{BLM96} G. R. Baker, X. Li, and A. C. Morlet. {\sl Analytic structure of two 1D-transport equations with nonlocal fluxes,} Phys. D {\bf 91} (1996) 349--375.

\bibitem{BKM} P. Biler, G. Karch,  R. Monneau. {\sl Nonlinear diffusion of dislocation density and self-similar solutions}. Comm. Math. Phys. {\bf 294} (2010), no. 1, 145--168. %MR2575479.

\bibitem{CChV} L.~A. Caffarelli, Ch.-H. Chan, and A. Vasseur.
{\sl  Regularity theory for nonlinear integral operators}, J. Amer. Math. Soc {\bf 24} (2011), 849--869.

\bibitem{CS07} L.~A. Caffarelli, L. Silvestre.
{ \sl An extension problem related to the fractional Laplacian.}
{Comm. Partial Differential Equations}  32 (2007),
no.~7-9, 1245--1260.

\bibitem{CSV} L.~A. Caffarelli. F. Soria, and J.~L. V\'azquez.
Regularity of solutions of the fractional  porous medium flow.
{\bf 15}, 5 (2013), 1701--1746. ArXiv 1201.6048v1 [math.AP].


\bibitem{CaVa09} L.~A. Caffarelli and J. L. V{\'a}zquez. {\sl Nonlinear porous medium flow
with fractional potential pressure},  Arch. Rational Mech. Anal. {\bf 202} (2011),
537--565.

\bibitem{CaVa11}  {L.~A. Caffarelli, and J. L. V{\'a}zquez}.
{\sl  Asymptotic behaviour of a porous medium equation
with fractional diffusion.}  DCDS-A {\bf 29}, no. 4 (2011), 1393--1404; A special
issue ``Trends and Developments in DE/Dynamics, Part III''.


\bibitem{Caffarelli-Vasseur} {L.~A. Caffarelli, A. Vasseur.}
{\sl Drift diffusion equations with fractional diffusion and the quasi-geostrophic equation.}
    Ann. of Math. (2) 171 (2010), no.~3, 1903--1930.


\bibitem{CFP2012} J. A. Carrillo, L. Ferreira, and J. C. Precioso.
{\sl A mass-transportation approach to a one dimensional fluid mechanics model with nonlocal velocity.} Adv. Math. {\bf 231} (2012), no. 1, 306–327.

\bibitem{CC08} A. Castro, D. C\'ordoba. {\sl Global existence, singularities and ill-posedness for a nonlocal flux.} Adv. Math. {\bf 219} (2008), no. 6, 1916--1936.


\bibitem{CCCF} D. Chae, A. C\'ordoba, D. C\'ordoba, M. A. Fontelos. {\sl Finite time singularities in a 1D model of the quasi-geostrophic equation.} Adv. Math. {\bf 194} (2005) 203--223.

\bibitem{dtdcnm04}
{ J.~Deslippe, R.~Tesdtrom, M.~S. Daw, D.~Chrzan, T.~Neeraj, and M.~Mills},
  {\em Dynamics scaling in a simple one-dimensional model of dislocation
  activity}, Phil. Mag., 84 (2004), pp.~2445--2454.

\bibitem{Head72} A. K. Head. {\sl Dislocation group dynamics I. Similarity solutions od the n-body problem}. Phil. Mag. {\bf 26} (1972), 43--53.

\bibitem{Landkof}  N. S. Landkof. {\sl ``Foundations of Modern Potential
Theory''.} Die Grundlehren der mathematischen Wissenschaften, {\bf 180}. Translated from the Russian by A. P. Doohovskoy. Springer, New York, 1972.

\bibitem{PQRV1} {\rm A.~de Pablo, F. Quir\'os, A. Rodriguez, J. L. V\'azquez. }\textit{A fractional porous medium equation} Adv. Math. \textbf{226} (2011), no. 2, 1378–1409.

\bibitem{PQRV2}{\rm A.~de Pablo, F. Quir\'os, A. Rodriguez, J. L. V\'azquez. } \textit{A general fractional porous medium equation},  Comm. Pure Applied Math. {\bf 65} (2012), 1242--1284.
\bibitem{Stein70} E. Stein. {\sl  ``Singular Integrals and Differentiability Properties of
Functions''}, Princeton University Press, Princeton, 1970.

\bibitem{SerVaz}{\rm S. Serfaty, J. L. V\'azquez,}
{\sl A Mean Field Equation as Limit of Nonlinear  Diffusion with Fractional Laplacian Operators},
Calc. Var. PDEs, to appear. ArXiv : 1205.632229 [math.AP].

\bibitem{Valdinoci} E. Valdinoci. {\it From the long jump random walk to the fractional Laplacian}.
Bol. Soc. Esp. Mat. Apl. {\bf 49} (2009), 33--44.

\bibitem{JLVAbel}{J.~L. V{\'a}zquez}.
{\sl Nonlinear Diffusion with Fractional Laplacian Operators.}
in ``Nonlinear partial differential equations: the Abel Symposium 2010'',
Holden, Helge  \& Karlsen, Kenneth H. eds., Springer, 2012. Pp. 271--298.

\bibitem{JLVsurvey2} {J. L. V{\'a}zquez}. {\sl Recent progress in the theory  of Nonlinear Diffusion with  Fractional Laplacian Operators}. In ``Nonlinear elliptic and parabolic differential equations'', Disc. Cont. Dyn. Syst. - S \ {\bf 7,} no. 4 (2014), 857--885. ArXiv:1401.3640.

\end{thebibliography}

\medskip

%\newpage

\noindent{\sc Addresses of the authors:}

\medskip

\noindent{\sc Luis A. Caffarelli}\newline
School of Mathematics, Univ. of Texas at Austin,\\
1 University Station, C1200, Austin, Texas 78712-1082. \newline
Second affiliation: Institute for Computational Engineering and Sciences.\newline
e-mail: {\tt  caffarel@math.utexas.edu}

\medskip

\noindent{\sc Juan Luis V{\'a}zquez}\newline
Departamento de Matem\'{a}ticas, Universidad Aut\'{o}noma de Madrid, \\
28049 Madrid, Spain. \ e-mail: {\tt juanluis.vazquez@uam.es}

\end{document}